\documentclass[preprint]{elsarticle}

\usepackage{amssymb}
\usepackage{amsmath}
\usepackage{amsthm}
\usepackage{enumerate}
\usepackage{graphicx}
\usepackage{amsfonts}

\journal{...}

\newtheorem{theorem}{Theorem}[section]
\newtheorem*{theoremA}{Theorem A}

\numberwithin{equation}{section}

\theoremstyle{definition}

\begin{document}

	\begin{frontmatter}

		\title{Overiteration of $d$-variate tensor product Bernstein operators: a quantitative result
		}

	 \author[1]{ Ana-Maria Acu}
	\author[2]{Heiner Gonska}

	\address[1]{Lucian Blaga University of Sibiu, Department of Mathematics and Informatics, Str. Dr. I. Ratiu, No.5-7, RO-550012  Sibiu, Romania,
		e-mail: anamaria.acu@ulbsibiu.ro; }
	\address[2]{University of Duisburg-Essen, Faculty of Mathematics, Bismarckstr. 90, D-47057 Duisburg, Germany, e-mail: heiner.gonska@uni-due.de and gonska.sibiu@gmail.com }

		\begin{abstract} 	
			{ Extending an earlier estimate for the degree of approximation of overiterated univariate Bernstein operators towards the same operator of degree one, it is shown that an analogous result holds in the $d$-variate case. The method employed can be carried over to many other cases and is not restricted to Bernstein-type or similar methods.
			} 
		\end{abstract}
	
		
	\end{frontmatter}

	\section{Introduction and historical remarks}
	
	The question behind this note is well-known. What is a classical Bernstein operator doing if its powers are raised to infinity?
	
	For the univariate version of this operator the answer is known. Already in 1966 the Dutch mathematician P.C. Sikkema proved  in the Romanian journal Mathematica (Cluj) that for  each function $f\in {\mathbb R}^{[0,1]}$ the powers $B_n^kf$, $n$ fixed, $k\to \infty$ converge to the linear function interpolating $f$ at $0$ and $1$ (see \cite {16}). Later on his result become known as the Kelisky-Rivlin (1967) or Karlin-Ziegler (1970) theorem (cf. \cite {11, 12}).
	
	However, even earlier T. Popoviciu \cite{Pop55} posed this problem in an (informal) problem book of 1955. We learned this from the note \cite{Alb79} by Albu cited by Precup \cite{Pre23}. The latter author also deals with multivariate operators but from a different point of view.
	
	Some notation is needed here. For $x\in[0,1]$, $n\geq 1$, and $f\in {\mathbb R}^{[0,1]}$ the Bernstein operator is given by
	\begin{align*}
		B_n(f,x)&:=\sum_{k=0}^nf\left(\dfrac{k}{n}\right)p_{n,k}(x)\\
		&:=\displaystyle\sum_{k=0}^{n}f\left(\dfrac{k}{n}\right){n\choose k}x^k(1-x)^{n-k}.
	\end{align*}

Thus $B_n$ is a polynomial operator, is linear and positive, reproduces all affine linear functions $l(x)=ax+b$, and for each $f$ the polynomial $B_nf$ is of degree $\leq n$.

Moreover, for any $k,n \in \mathbb N$, Gonska et al. \cite{8} proved in 2006, extending earlier work of Nagel \cite{14} and Gonska \cite{6},

\begin{align}
	|B_n^k(f,x)-B_1(f,x)
	&\leq \dfrac{9}{2}\omega_2\left(f;\sqrt{x(1-x)\left(1-\dfrac{1}{n}\right)^k}\right), x\in [0,1].\label{eq.*}
\end{align}
Here $\omega_2(f;\cdot)$ is the classical second order modulus of $f$. Hence the right hand side converges to $0$ as $n$ is fixed and $k\to\infty$ (some more general situations are possible). It also shows that the powers are interpolatory at $0$ and $1$ and keep reproducing linear functions. Moreover, the convergence is uniform with respect to $\|\cdot\|_{\infty}$.

When it comes to multivariate Bernstein operators, all the time operators on generalized simplices or hypercubes are meant. While for simplices the convergence of powers was investigated by, e.g., Wenz \cite{17} and many others, the hypercube case remained allegedly open until a 2009 article of Jachymski \cite{10} appeared. However, for the bivariate case a paper by Agratini and Rus was published already in 2003, cf. \cite{AgRu2003b}. 

In this note we will use the term tensor product although in other publications one might find 'product of parametric extensions' meaning exactly the same (see, e.g., \cite{7}).

Using functional-analytic methods Jachymski showed the following. For $l,m\geq 1$ let the bivariate tensor product operator 
\begin{align*}
	\left((B_l\otimes  B_m)f\right)(x,y)&:=
\left(_sB_l\circ\,\! _tB_m\right)\left(f(s,t);x,y\right)
	\end{align*}
be given by

$$\displaystyle\sum_{i=0}^l\sum_{j=0}^m f\left(\dfrac{i}{l},\dfrac{j}{m}\right)p_{l,i}(x)\cdot p_{m,j}(y),\,\, f\in C([0,1]^2),\,x,y\in[0,1].  $$

\begin{theoremA} For any $l,m\in{\mathbb N}$ fixed, the sequence $\left((B_l \otimes B_m)^n\right)_{n\in{\mathbb N}}$ uniformly converges to the operator $L$ (independent of $l$ and $m$) given by the following formula for $f\in C([0,1]^2)$ and $x,y\in [0,1]$:
	\begin{align*}
		(Lf)(x,y)\\
		&=f(0,0)+(f(1,0)-f(0,0))x+(f(0,1)-f(0,0))y+(f(0,0)+f(1,1)-f(1,0)-f(0,1))xy\\
		&=(1-x,x)\left(\begin{array}{cc}
			f(0,0) &f(0,1)\\
			f(1,0) & f(1,1)
			\end{array} \right){1-y\choose y}.
	\end{align*}
In other words, $Lf=(B_1 \otimes B_1)f$.

\end{theoremA}

Jachymski \cite{10} also gave the limit of $n$-powers of $d$-variate Bernstein operators, i.e., of
$$\left((B_{l_1} \otimes \dots \otimes B_{l_d})f\right)(x_1,\dots,x_d)  $$
$$ =\left(_{s_1}B_{l_1}\circ\dots\circ\,\!_{s_d}B_{l_d}\right)\left(f(s_1,\dots,s_d);x_1,\dots,x_d\right).  $$
They map $C([0,1]^d)$ into $\Pi_{l_1,\dots,l_d}$, the space of $d$-variate polynomials of total degree $\leq \displaystyle\sum_{\delta=1}^dl_{\delta}$.

The limiting operator in this case is 
$$ (Lf)(x_1,\dots,x_d)=\displaystyle\sum_{(\epsilon_1,\dots,\epsilon_d)\in V}f(\epsilon_1,\dots,\epsilon_d)p_{\epsilon_1}(x_1)\cdot\dots\cdot p_{\epsilon_d}(x_d), $$
where $V=\{0,1\}^{\{1,\dots,d\}}$, and for $s\in [0,1]$, $p_0(s):=1-s$ and $p_1(s):=s$. Thus $L$ equals $B_1\otimes...\otimes B_1.$

In the present note we will show first that the fixpoint approach of (Agratini and) Rus also works in the $d$-variate case. Our main emphasis is on the quantitative situation where we will demonstrate how the pointwise $\omega_2$-result may be carried over to $d$ dimensions.

\section{The non-quantitative approach of Agratini and Rus revisited}
As mentioned above, Jachymski used a functional-analytic framework to derive his result. Here we show that a more elementary approach does the job as well. We recall the three papers by Rus and Agratini $\&$ Rus and present their approach for $d$ dimensions.

Some reminders concerning $d$-variate hypercubes are in order. More details are available in the German Wikipedia, keyword "Hyperwürfel" \cite{18}. Such a hypercube in $d$ dimensions possesses ${d\choose 0}2^{d-0}=2^d$ 0-dimensional boundary elements (vertices), in the bivariate case these are the 4 corners of $[0,1]^2$. Adopting the above notation these are all $d$-tuples 
$$ (\epsilon_1,\dots,\epsilon_d)\in V,\,\, V=\{0,1\}^{\{1,\dots,d\}}. $$
We will now follow Rus' proof of his Theorem 1. First introduce the sets
$$ X_{\alpha_1,\dots,\alpha_d}=\{f\in C([0,1]^d)\,:\, f(\epsilon_1)=\alpha_1,\dots, f(\epsilon_d)=\alpha_d\},\,\, (\epsilon_1,\dots,\epsilon_d)\in V,\, \alpha_1,\dots,\alpha_d\in{\mathbb R}.  $$
Note that
\begin{itemize}
	\item[(a)] $X_{\alpha_1,\dots,\alpha_d}$ is a closed subset of $C([0,1]^d)$;
	\item[(b)] $X_{\alpha_1,\dots,\alpha_d}$ is an invariant subset of $B_{l_1}\otimes \dots \otimes B_{l_d}$, for all $\alpha_1,\dots,\alpha_d\in {\mathbb R}$ and $l_1,\dots,l_d \in {\mathbb N}$;
	\item[(c)] $C([0,1]^d)=\displaystyle\bigcup_{\alpha_1,\dots,\alpha_d\in {\mathbb R}}X_{\alpha_1,\dots,\alpha_d}$ is a partition of $C([0,1]^d)$.
\end{itemize}
Next it is shown that
$$ \left.(B_{l_1}\otimes\dots \otimes B_{l_d})\right|_{X_{\alpha_1,\dots, \alpha_d}} $$
maps $X_{\alpha_1,\dots, \alpha_d}$ onto itself and is a contraction.

For $f,g\in X_{\alpha_1,\dots, \alpha_d}$ we have
$$ \left| \left((B_{l_1}\otimes\dots \otimes B_{l_d})f\right)(x_1,\dots,x_d) -\left((B_{l_1}\otimes\dots \otimes B_{l_d})g\right)(x_1,\dots,x_d)\right|  $$
$$ =\left| \displaystyle\sum_{\lambda_1=0}^{l_1}\dots \sum_{\lambda_d=0}^{l_d}(f-g)\left(\dfrac{\lambda_1}{l_1},\dots,\dfrac{\lambda_d}{l_d}\right)p_{l_1,\lambda_1}(x_1)\cdot\dots\cdot p_{l_d,\lambda_d}(x_d)\right| $$
$$  \leq \sum_{(\lambda_1,\dots,\lambda_d)\in\{0,\dots,l_1\}\times\dots \times\{0,\dots,l_d\}\setminus V }\left|(f-g)\left(\dfrac{\lambda_1}{l_1},\dots,\dfrac{\lambda_d}{l_d}\right)p_{l_1,\lambda_1}(x_1)\cdot\dots\cdot p_{l_d,\lambda_d}(x_d)\right| $$
$$ \leq \|f-g\|_{\infty}\sum_{(\lambda_1,\dots,\lambda_d)\in\{0,\dots,l_1\}\times\dots \times\{0,\dots,l_d\}\setminus V }p_{l_1,\lambda_1}(x_1)\cdot\dots\cdot p_{l_d,\lambda_d}(x_d)$$
$$ \leq \|f-g\|_{\infty}\left(1-\min \sum_{(\lambda_1,\dots,\lambda_d)\in V }p_{l_1,\lambda_1}(x_1)\cdot\dots\cdot p_{l_d,\lambda_d}(x_d)\right)  $$
$$= \|f-g\|_{\infty}\cdot \left(1-\min\left\{\left[(1-x_1)^{l_1}+x_1^{l_1}\right]\cdot \dots \cdot \left[(1-x_d)^{l_d}+x_d^{l_d}\right]\right\}\right)  $$
$$\leq  \|f-g\|_{\infty}\cdot \left(1-\dfrac{1}{\prod_{\delta=1}^d2^{l_{\delta}-1}}\right)<1. $$
Thus $B_{l_1}\otimes\dots \otimes B_{l_d} $ on $X_{\alpha_1,\dots,\alpha_d}$ is a contraction for all $\alpha_1,\dots,\alpha_d\in{\mathbb R}$. On the other hand, $(Lf)(x_1,\dots,x_d)$ is a fixed point of $B_{l_1}\otimes\dots \otimes B_{l_d}$.

So $f\in C([0,1]^d)$ is in $X_{f(\epsilon_1),\dots,f(\epsilon_d)}$ and from the contraction principle we have
$$\displaystyle\lim_{n\to\infty}\left(B_{l_1}\otimes\dots \otimes B_{l_d}\right)^nf=Lf.  $$

We summarize our observation in
\begin{theorem}
	(Jachymski \cite{10})
For fixed $l_1,\dots,l_d\in {\mathbb N}=\{1,2,\dots\}$ one has
$$\displaystyle\lim_{n\to\infty}\left(B_{l_1}\otimes\dots \otimes B_{l_d}\right)^nf=Lf \textrm{ uniformly}.  $$
Here $B_{l_1}\otimes\dots \otimes B_{l_d}$ is the $d$-variate tensor product operator on $C([0,1]^d)$ and
$$ (Lf)(x_1,\dots,x_d)=\displaystyle\sum_{(\epsilon_1,\dots,\epsilon_d)}f(\epsilon_1,\dots,\epsilon_d)p_{\epsilon_1}(x_1)\cdot\dots\cdot p_{\epsilon_d}(x_d),\,\, V=\{0,1\}^{\{1,\dots,d\}}.  $$
\end{theorem}
In particular, for $d=2$ we have the representation of Theorem A.

\section{The Zhuk extension in the bi- and $d$-variate cases}

Since  the articles of Zhuk  \cite{19} and Gonska $\&$ Kovacheva \cite{9} are hard to obtain, we briefly describe the extension in the univariate situation, then carry it over to the bivariate case and finally show what has to be done in $d$ variables.

\subsection{Zhuk construction-univariate case}
For $f\in C[0,1]$ and $0<h\leq\dfrac{1}{2}(b-a)$ define $f_h:[a-h,b+h]\to {\mathbb R}$ by
$$ f_h(x):=\left\{\begin{array}{l}P_{-}(x),\,\, a-h\leq x<a,\\
\vspace{-0.3cm}\\
f(x),\,\, a\leq x\leq b,\\
\vspace{-0.3cm}\\
P_{+}(x),\,\, b<x\leq a+h.
\end{array}\right.  $$
\begin{align*}
&	\|f-P_{-}\|_{C[a,a+2h]}=E_1(f;a,a+2h),\\
&	\|f-P_{+}\|_{C[b-2h,b]}=E_1(f;b-2h,b).
\end{align*}
Here $P_{-}$ and $P_{+}$ denote the best approximations in $\Pi_1$ on the intervals indicated and with respect to the uniform norm.

Zhuk put
$$ S_{h}(f;x):=\dfrac{1}{h}\int_{-h}^h \left(1-\dfrac{|t|}{h}\right)f_h(x+t)dt,\, x\in[a,b].  $$

He showed \cite[Lemma 1]{19}:
For $f\in C[a,b], 0<h\leq \dfrac{1}{2}(b-a)$,
\begin{align*}
	& \|f-S_hf\|_{\infty}\leq \dfrac{3}{4}\omega_2(f;h),\\
	&\|(S_hf)^{\prime\prime}\|_{L_{\infty}[a,b]}\leq \dfrac{3}{2}h^{-2}\omega_2(f;h). 
\end{align*}

\subsection{Construction of the bivariate Zhuk extension}
Let $f\in C([0,1]^2)$. On a fixed $y$-level we extend the partial function $f_{y}(x) = f(\cdot ,y)$ from $[0,1]\times \{y\}$ to $[-h, 1+h]\times\{y\} $ in complete analogy to the univariate case. After integration, for each $y\in [0,1]$, we obtain
$$ S_{h}(f_y;x):=\dfrac{1}{h}\int_{-h}^h\left(1-\dfrac{|t|}{h}\right)(f_y)_h(x+t)dt,\,\, x\in[0,1],  $$
satisfying for $0<h\leq \dfrac{1}{2}$:
\begin{align*}
&	\|f_y-S_hf_y\|_\infty \leq \dfrac{3}{4}\omega_2(f_y;h),\\
& \|(S_hf_y)^{\prime\prime}\|_{L_{\infty}[0,1]}\leq \dfrac{3}{2}h^{-2}\omega_2(f_y;h).
\end{align*}
(On each $y$-level we could have even chosen $h_y$ with $0<h_y\leq \dfrac{1}{2}$).

The same procedure we carry out for $f_x(y), y\in [0,1]$, producing functions $S_hf_x$ such that

\begin{align*}
	&	\|f_x-S_hf_x\|_C\leq \dfrac{3}{4}\omega_2(f_x;h),\\
	& \|(S_hf_x)^{\prime\prime}\|_{L_{\infty}[0,1]}\leq \dfrac{3}{2}h^{-2}\omega_2(f_x;h).
\end{align*}
This can be done for all $x\in [0,1]$.

More explicitly,
\begin{align*}
	\omega_2(f_y;h)&=\sup\{|f_y(x-\delta)-2f_y(x)+f_y(x+\delta)|:\,|\delta|\leq h, x\pm \delta\in[0,1]\}\\
	&=\sup\{|f(x-\delta,y)-2f(x,y)+f(x+\delta,y)|:\,|\delta|\leq h, x\pm \delta\in[0,1]\}\\
	&\leq \sup_{y\in[0,1]}\sup\{ |f(x-\delta,y)-2f(x,y)+f(x+\delta,y)|:\,|\delta|\leq h, x\pm \delta\in[0,1]\}\\
	&=\omega_2(f;h,0).
\end{align*}
Also, $\omega_2(f_x;h)\leq \omega_2(f;0,h)$.

The quantities $\omega_2(f;h,0)$ and $\omega_2(f;0,h)$ are called "partial moduli of smoothness". We have thus constructed auxiliary extensions of $f_y(\cdot)$, $y\in[0,1]$, and $f_x(\ast)$, $x\in[0,1]$, on the domain shown below

	\begin{figure}[htbp]
	\centering
	\includegraphics[height=60mm,keepaspectratio]{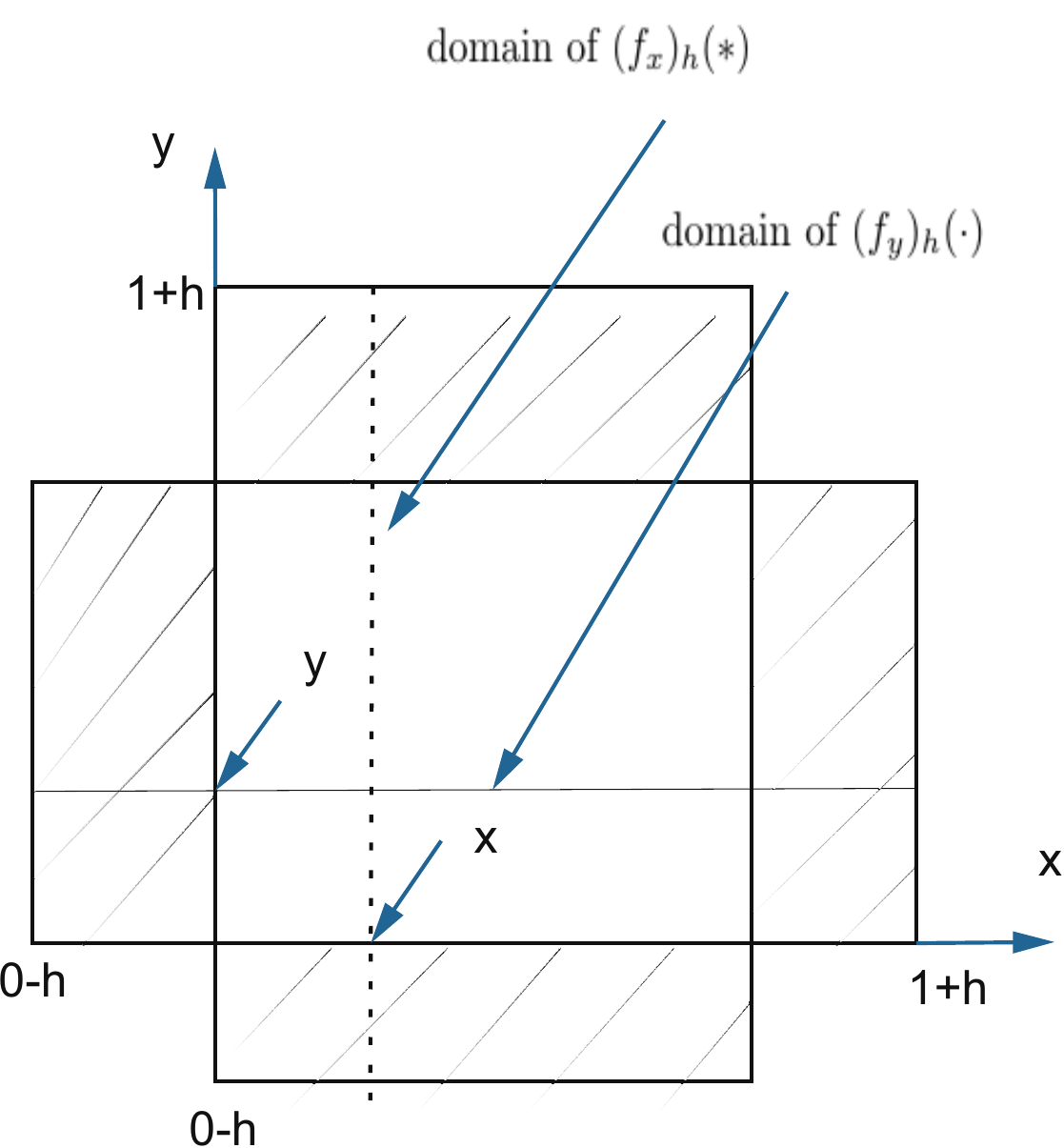}
	\caption{ } 
	\label{fig:4}
\end{figure}

$S_h(f_y;\cdot)$  and $S_h(f_x;\ast)$ are given on the inner (white) square only.

\subsection{Zhuk extension, $d$-variate case}

The construction described for the bivariate case can be easily generalized for $d\geq 3$ dimensions. To this end fix $d-1\geq 2$ variables, say $s_2,\dots, s_d$. Then extend the partial function $f_{s_2,\dots,s_d}(s_1)$, $0\leq s_1\leq 1$, to $-h\leq s_1\leq 1+h$, $0<h\leq \dfrac{1}{2}$, and define
$$ S_h(f_{s_2,\dots,s_d})(s_1):=\dfrac{1}{h}\int_{-h}^h\left(1-\dfrac{|t|}{h}\right)\cdot\left(f_{s_2,\dots,s_d}\right)_h(s_1+t)dt.  $$
This gives
$$\|f_{s_2,\dots,s_d}-S_hf_{s_2,\dots,s_d}\|_{\infty}\leq \dfrac{3}{4}\omega_2(f_{s_2,\dots,s_d};h),  $$
$$ \|(S_hf_{s_2,\dots,s_d})^{\prime\prime}\|_{L_{\infty}[0,1]}\leq \dfrac{3}{2}h^{-2}\omega_2(f_{s_2,\dots,s_d};h),  $$
for each fixed $s_2,\dots,s_d\in[0,1]$.
Moreover, a common upper bound is
$$ \omega_2(f_{s_2},\dots,s_d;h)\leq \omega_2(f;h,0,\dots,0),\textrm{ for all } s_2,\dots,s_d\in[0,1],  $$
and a corresponding inequality holds for any other choice of $s_{\delta}$, $2\leq \delta\leq d$.

\section{An estimate for $d$-variate tensor product Bernstein operators}

We first recall our 2006 estimate for the univariate case:
$$ |B_l^n(f;x)-B_1(f;x)|\leq\dfrac{9}{4}\omega_2\left(f;\sqrt{x(1-x)\left(1-\dfrac{1}{l}\right)^n}\right). $$
In two dimensions, it can be easily derived that
$$\left|(B_l\otimes B_m)^n (f;x,y)-(B_1\otimes B_1) (f;x,y)\right|$$
$$ \left|\left[(B_l^n-B_1)\otimes (B_m^n-B_1)\right](f;x,y) \right| $$
$$ \leq \dfrac{9}{4}\left[\omega_2\left(f;\sqrt{x(1-x)\left(1-\dfrac{1}{l}\right)^n},0\right)+\omega_2\left(f;0,\sqrt{y(1-y)\left(1-\dfrac{1}{m}\right)^n}\right)\right].  $$
This extends to $d$ dimensions. Here we have 
$$ \left| \left(_{s_1}B_{l_1}\circ\dots\circ\,\!_{s_d}B_{l_d}\right)^n\left(f(s_1,\dots,s_d);x_1,\dots,x_d\right) 
-\left(_{s_1}B_{1}\circ\dots\circ\,\!_{s_d}B_{1}\right)\left(f(s_1,\dots,s_d);x_1,\dots,x_d\right)\right| $$
$$ \leq \dfrac{9}{4}\sum_{\delta=1}^{d} 
\omega_2\left(f;0,\dots ,0,\sqrt{x_{\delta}(1-x_{\delta})\left(1-\dfrac{1}{l_{\delta}}\right)^n},0, \dots,0\right). $$
	For $d$ dimensions it is, without additional effort, possible to show 

$$ \left| \left(_{s_1}B_{l_1}^{n_1}\circ\dots\circ\,\!_{s_d}B_{l_d}^{n_d}\right)\left(f(s_1,\dots,s_d);x_1,\dots,x_d\right) 
-\left(_{s_1}B_{1}^{n_1}\circ\dots\circ\,\!_{s_d}B_{1}^{n_d}\right)\left(f(s_1,\dots,s_d);x_1,\dots,x_d\right)\right| $$
$$ =\left[(_{s_1}B_{l_1}^{n_1}-_{s_1}B_{1})\circ\dots\circ (_{s_d}B_{l_d}^{n_d}-_{s_d}B_{1})\right]\left(f(s_1,\dots,s_d);x_1,\dots,x_d\right) $$
$$ \leq \dfrac{9}{4}\sum_{\delta=1}^{d} 
\omega_2\left(f;0,\dots 0,\sqrt{x_{\delta}(1-x_{\delta})\left(1-\dfrac{1}{l_{\delta}}\right)^{n_{\delta}}},0,\dots ,0\right). $$
Note that for $n=n_1=\dots=n_{\delta}$ the difference from above becomes 
$$\left(_{s_1}B_{l_1}\otimes\dots\otimes\,\!_{s_d}B_{l_d}\right)^n
-\left(_{s_1}B_{1}\otimes\dots\otimes\,\!_{s_d}B_{1}\right)$$
and is this the multivariate quantity considered by Jachymski. However, there is no need to restrict oneself to this case.

\vspace{.4cm}

\section{Optimality}

Questions are in order in how far our estimates are "optimal".

1. The constant $\dfrac{9}{4}$ appearing repeatedly in this note most likely is not. There is need for work in this direction.

2. If the function $f$ is $d$-linear, then the sum of $d$ $ \omega_2$ -terms equals zero.
If the sum is zero, then each of its terms does so. This may occur if

(i) $(x_1, \dots, x_d)$ is at a 'corner' of the hypercube, and/or 

(ii) $l_\delta$, the degree of $B_{l_\delta}$, is equal to 1 for $1 \leq \delta \leq d$.

In any other case f must be $d$-linear to fulfill the condition $\omega_2 (f;\dots)=0$ for all $d$ terms and for an interior point of the hypercube while $l_{\delta} \ge 2, 1 \leq \delta \leq d.$

From (i) and (ii) it is evident that the sum of $d$ $\omega_2$-terms is the correct expression for tensor product Bernstein approximation over a (generalized) hypercube.

\section{Concluding remark}

It should have become clear that our, or a similar approach, may be used to prove analogous results for many other operator sequences  (which different authors may consider). We feel that sums of partial moduli of smoothness are among the right tools for tensor product approximation since they show the mutual independence of the variables. Nonetheless, even better pointwise results are available but do not really contribute to a better understanding.


\begin{thebibliography}{99}
	
	
	
	
\bibitem{3}	O. Agratini, I. Rus, 
	Iterates of a class of discrete linear operators via contraction principle.
	Commentat. Math. Univ. Carol. 44, no. 3, 555-563 (2003). 
	
\bibitem{AgRu2003b} O. Agratini, I. Rus, 
Iterates of some bivariate approximation process via weakly Picard operators.
Nonlinear Anal. Forum 8, no. 2 (2003), 159-168.

\bibitem{Alb79} M. Albu,
Asupra convergen\c tei iteratelor unor operatori liniari \c si m\u argini\c ti.
Sem. Itin. Ec. Fun. Approx. Convex., Cluj-Napoca 1979, 6-16.
	
	
	
\bibitem{6}	H. Gonska,
	Quantitative Aussagen zur Approximation durch positive lineare Operatoren, Dissertation, Universität Duisburg 1979.
	
\bibitem{7}	H.H. Gonska,
	Products of parametric extensions: refined estimates.
	In: Proc. 2nd Int. Conf. on “Symmetry and Antisymmetry” in Mathematics, Formal Languages and Computer Science (ed. by G.V. Orman and D. Bocu), 1-15. Brasov: Editura Universit\u a\c tii “Transilvania” 2000. 
	
\bibitem{8}	H. Gonska, D. Kacsó, P. Pitul,
	The degree of convergence of over-iterated positive linear operators,
	J. of Applied Functional Analysis 1 (2006), 403-423.
	
\bibitem{9}	H. Gonska, R. Kovacheva,
	The second order modulus revisited: remarks, applications, problems,
	Confer. Sem. Mat. Univ. Bari 257 (1994), 1-32.
	
\bibitem{10}	J. Jachymski, 
	Convergence of iterates of linear operators and the Kelisky-Rivlin type theorems,
	Studia Mathematica 195, no. 2 (2009), 99-112.

\bibitem{12}	S. Karlin, Z. Ziegler,
	Iteration of positive approximation operators,
	J. Approx. Th. 3 (1970), 310-339.
		
\bibitem{11}	R.P. Kelisky, T.J. Rivlin,
	Iterates of Bernstein polynomials,
	Pac. J. Math. 21 (1967), 511-520.
	
	
\bibitem{14}	J. Nagel, 
	Sätze Korovkinschen Typs für die Approximation linearer positiver Operatoren,
	Dissertation, Universität Essen 1978. 
	
\bibitem{Pop55} T. Popoviciu,
Problem posed on Dec. 6, 1955. In:	Caietul de probleme al Catedrei de analiz\u a, n. 58.

\bibitem{Pre23} R. Precup,
On the iterates of uni- and multidimensional operators.
Bull. Transilv. Univ. Brasov, Ser. III: Math. and Comp. Sci. 3 (65), No. 2 (2023),143-152.
	
\bibitem{15}	I.A. Rus,
	Iterates of Bernstein operators, via contraction principle,
	J. Math. Anal. Appl. 292 (2004), 259-261.
	
\bibitem{16}	P.C. Sikkema,
	Über Potenzen von verallgemeinerten Bernstein-Operatoren, 
	Mathematica, Cluj 8(31), 173-180 (1966).
	
\bibitem{17} H.-J. Wenz, 
	On the limits of (linear combinations of) iterates of linear operators.
	J. Approx. Th. 89, no. 2 (1997), 219-237.
	
\bibitem{18}	Wikipedia (German),
  https://de.wikipedia.org/wiki/Hyperw\"urfel (Jan. 19, 2024).
	
\bibitem{19}	V. V. Zhuk,
	Functions of the Lip 1 class and S. N. Bernstein’s polynomials (Russian),
	Vestn. Leningr. Univ., Math. 22, no. 1 (1989), 38-44; translation of Vestn. Leningr. Univ., Ser. I 1989, no. 1, 25-30 (1989).
	

	
	
\end{thebibliography}
\end{document}